\documentclass[a4paper]{amsart}

\usepackage[english]{babel}
\usepackage[fixlanguage]{babelbib}
\usepackage[utf8]{inputenc}
\usepackage{amsfonts,amssymb,amsmath,amsthm}

\theoremstyle{definition}
\newtheorem{defi}{Definition}[section]

\theoremstyle{plain}
\newtheorem{thm}[defi]{Theorem}
\newtheorem*{thm*}{Theorem}
\newtheorem{lemma}[defi]{Lemma}

\theoremstyle{remark}

\newcommand{\Fp}{\mathbb F_p}
\newcommand{\Z}{ \mathbb{Z} }
\newcommand{\N}{ \mathbb{N} }

\newcommand{\PP}{ \mathbb{P} }

\newcommand{\dd}{ \mathsf{d} }
\newcommand{\C}{ \mathsf{C} }
\newcommand{\lbar}{ \overline }
\providecommand{\length}[1]{\lvert#1\rvert}
\providecommand{\card}[1]{\lvert#1\rvert}
\providecommand{\lrcard}[1]{\Big\lvert#1\Big\rvert}
\providecommand{\val}{ \mathsf{v} }

\newcommand{\be}{\begin{equation}}
\newcommand{\ee}{\end{equation}}

\newcommand{\ber}{\begin{eqnarray}}
\newcommand{\eer}{\end{eqnarray}}
\newcommand{\nn}{\nonumber}

\DeclareMathOperator{\ord}{ord}

\renewcommand{\t}{\, | \,}

\begin{document}

\address{Institut f\"ur Mathematik und Wissenschaftliches Rechnen \\
Karl-Franzens-Uni\-ver\-si\-t\"at Graz \\
Heinrichstra\ss e 36\\
8010 Graz, Austria} \email{06smertn@edu.uni-graz.at}

\title{On the Davenport constant and group algebras}

\author{Daniel Smertnig}


\thanks{This work was supported by the Austrian Science Fund FWF,
Doctoral Program C75-N13 Discrete Mathematics}.


\begin{abstract}
For a finite abelian group $G$ and a splitting field $K$ of $G$, let
$\mathsf d (G, K)$ denote the largest integer $l \in \N$ for which
there is a sequence $S = g_1 \cdot \ldots \cdot g_l$ over $G$ such
that $(X^{g_1} - a_1) \cdot \ldots \cdot (X^{g_l} - a_l)  \ne  0 \in
K[G]$ for all $a_1, \ldots, a_l \in K^{\times}$. If $\mathsf D (G)$
denotes the Davenport constant of $G$, then there is the
straightforward inequality $\mathsf D (G)-1 \le \mathsf d (G, K)$.
Equality holds for a variety of groups, and a standing conjecture of
W. Gao et.al. states that equality holds for all groups. We offer
further groups for which equality holds, but we also give the first
examples of groups $G$ for which $\mathsf D (G) -1 < \mathsf d (G,
K)$ holds. Thus we disprove the conjecture.
\end{abstract}

\maketitle

\bigskip
\section{Introduction and Main Result}
\bigskip

Let $G$ be an additive finite abelian group. For a (multiplicatively
written) sequence $S = g_1 \cdot \ldots \cdot g_l$ over $G$, $|S| =
l$ is called the length of $S$, and $S$ is said to be zero-sum free
if $\sum_{i \in I}g_i \ne 0$ for every nonempty subset $I \subset
[1, l]$. Let $\mathsf d (G)$ denote the maximal length of a zero-sum
free sequence over $G$. Then $\mathsf d (G) + 1$ is the Davenport
constant of $G$, a classical constant from Combinatorial Number
Theory (for surveys and historical comments, the reader is referred
to \cite{Ga-Ge06b}, \cite[Chapter 5]{Ge-HK06a}, \cite{Ge09a}). In
general, the precise value of  $\mathsf d (G)$ (in terms of the
group invariants of $G$) and the structure of the extremal sequences
is unknown, see \cite{Le-Sc07a, Bh-SP07a, Ra-Sr-Th08, Gi09b, Gi09c,
Ga-Ge-Gr09a, Sc09e, Su09a, Ge-Li-Ph10a} for recent progress.

\smallskip
Group algebras $R [G]$ - over suitable commutative rings $R$ - have
turned out to be powerful tools for a great variety of questions
from combinatorics and number theory, among them the Davenport
constant. We recall the definition of an invariant (involving group
algebras) which was used for the investigation of the Davenport
constant since the 1960s.

\smallskip
For a commutative ring $R$, let \ $\mathsf d (G, R) \in \mathbb N
\cup \{\infty\}$ \ denote the supremum of all $l \in \mathbb N$
having the following property:

\begin{enumerate}
\item[]
There is some sequence \ $S = g_1 \cdot \ldots \cdot g_l$ \ of
length $l$ over $G$ such that
\[
(X^{g_1} - a_1) \cdot \ldots \cdot (X^{g_l} - a_l)  \ne  0 \in R[G]
\quad \text{for all} \quad a_1, \ldots , a_l \in R \setminus \{0\}
\,.
\]
\end{enumerate}

\smallskip

\noindent If $S$ is zero-sum free, $R$ is an integral domain, $a_1,
\ldots, a_l \in R \setminus \{0\}$ and
\[
f \ = \ (X^{g_1} - a_1) \cdot \ldots \cdot (X^{g_l} - a_l) \ = \
\sum_{g \in G} c_g X^g \ \,,
\]
then $c_0 \ne 0$. Hence $f \ne 0$, and it follows that
\[
\mathsf d (G) \le \mathsf d (G, R)\,.
\]
The following \,Theorem {\bf A}\, was proved by P. van Emde Boas,
D. Kruyswijk and J.E. Olson in the 1960s (in fact, they did not
explicitly define the invariants $\mathsf d (G, K)$ but got these
results implicitly). Historical remarks and proofs in the present
terminology may be found in \cite[Section 2.2]{Ge09a} and
\cite[Theorem 5.5.9]{Ge-HK06a}; see also \cite{Ga-Ge-HK09}.

\medskip

\noindent
{\bf Theorem A.} {\it Let $G$ be a finite abelian group
with $\exp (G) = n \ge 2$.
\begin{enumerate}
\item Let $K$ be a splitting field of $G$ with $\text{\rm char} (K) \nmid \exp (G)$.
      Then
      \[
      \mathsf d (G, K) \le (n-1) + n \log \frac{|G|}{n} \,.
      \]

\smallskip
\item
If $G$ is a $p$-group, then $\mathsf d (G) = \mathsf d (G,
      \mathbb Z/ p \mathbb Z)$.
\end{enumerate}
}

\smallskip
Note that for a cyclic group $G$ of order $n$,  the above upper
bound implies that $\mathsf d (G) = \mathsf d (G, K) = n-1$, since $\mathsf d(\C_n) \ge n-1$ can easily be seen. Only
recently, W. Gao and Y. Li  showed that $\mathsf d (\C_2 \oplus
\C_{2n}) = \mathsf d (\C_2 \oplus \C_{2n}, K)$ (\cite[Theorem
3.3]{Ga-Li09a}). We  extend their result, but we also show that
Conjecture 3.4 in \cite{Ga-Li09a}, stating that $\mathsf d (G) =
\mathsf d (G, K)$  for all groups $G$, does not hold. Here is the
main result of the present paper.

\medskip
\begin{thm} \label{thm:main}
  Let $G = \C_p \oplus \C_{pn}$ with $p \in \PP$, $n \in \N$ and let $K$ be a splitting field of $G$.
  \begin{enumerate}
    \smallskip
    \item If $p \le 3$, then $\mathsf d (G) = \dd(G,K)$.

    \smallskip
    \item If $p \ge 5$ and $n \ge 2$, then $\mathsf d (G) < \dd(G,K)$.
  \end{enumerate}
\end{thm}

\bigskip
\section{Preliminaries}
\bigskip

Let $\mathbb N$ denote the set of positive integers, $\PP \subset
\N$ the set of prime numbers, and let $\mathbb N_0 = \mathbb N \cup
\{ 0 \}$. For real numbers $a, b \in \mathbb R$, we set $[a, b] = \{
x \in \mathbb Z \mid a \le x \le b\}$.  For $n \in \mathbb N$ and $p
\in \PP$, let $\C_n$ denote a cyclic group with $n$ elements,
$\mathsf v_p (n) \in \N_0$ the $p$-adic valuation of $n$ with
$\mathsf v_p (p) = 1$ and $\Fp = \Z / p \Z$ the finite field with
$p$ elements.

\smallskip
Let $G$ be an additive finite abelian group.  Suppose that  $G \cong
C_{n_1} \oplus \ldots \oplus C_{n_r}$ with $1 < n_1 \t \ldots \t
n_r$.  Then $r = \mathsf r (G)$ is the {\it rank} of $G$, $n_r =
\exp (G)$ is the {\it exponent} of $G$, and we define $\mathsf d^*
(G) = \sum_{i=1}^r (n_i-1)$. If $|G| = 1$, then the exponent $\exp
(G) = 1$, the rank $\mathsf r (G) = 0$,
and we set $\mathsf d^* (G) = 0$.
If $A, B \subset G$ are nonempty subsets, then $A + B = \{ a + b
\mid a \in A, b \in B \}$ is their sumset. We will make use of a
Theorem of Cauchy-Davenport which runs as follows (for a proof see
\cite[Cor. 5.2.8.1]{Ge-HK06a}).

\medskip
\begin{lemma} \label{2.1}
  Let $G$ be a cyclic group of order $p \in \PP$ and let $A,B \subset G$ be nonempty subsets.
  Then $\card{A+B} \ge \min\{ \card{A}+\card{B}-1, p \}$.
\end{lemma}

\medskip
{\bf Sequences over groups.} Let $\mathcal F(G)$ be the
(multiplicatively written)  free abelian monoid with basis $G$. The
elements of $\mathcal F(G)$ are called \ {\it sequences} \ over $G$.
We write sequences $S \in \mathcal F (G)$ in the form
\[
S =  \prod_{g \in G} g^{\mathsf v_g (S)}\,, \quad \text{with} \quad
\mathsf v_g (S) \in \mathbb N_0 \quad \text{for all} \quad g \in G
\,.
\]
We call \ $\mathsf v_g (S)$  the \ {\it multiplicity} \ of $g$ in
$S$, and we say that $S$ \ {\it contains} \ $g$ \ if \ $\mathsf v_g
(S) > 0$.   A sequence $S_1 $ is called a \ {\it subsequence} \ of
$S$ \ if \ $S_1 \, | \, S$ \ in $\mathcal F (G)$ \ (equivalently, \
$\mathsf v_g (S_1) \le \mathsf v_g (S)$ \ for all $g \in G$). If a
sequence $S \in \mathcal F(G)$ is written in the form $S = g_1 \cdot
\ldots \cdot g_l$, we tacitly assume that $l \in \mathbb N_0$ and
$g_1, \ldots, g_l \in G$. For a sequence
\[
S \ = \ g_1 \cdot \ldots \cdot g_l \ = \  \prod_{g \in G} g^{\mathsf
v_g (S)} \ \in \mathcal F(G) \,,
\]
we call
\[
|S| = l = \sum_{g \in G} \mathsf v_g (S) \in \mathbb N_0 \qquad
\text{the \ {\it length} \ of \ $S$ \ and}
\]
\[
\sigma (S) = \sum_{i = 1}^l g_i = \sum_{g \in G} \mathsf v_g (S) g
\in G \qquad \text{the \ {\it sum} \ of \ $S$}\,.
\]
The sequence \ $S$ \ is called a {\it zero-sum sequence} if $\sigma
(S) = 0$, and it is called {\it zero-sum free} if $\sum_{i \in I}
g_i \ne 0$ for all $\emptyset \ne I \subset [1, l]$ (equivalently,
if there is no nontrivial zero-sum subsequence). We denote by
\begin{itemize}
\item $\mathsf D (G)$ the smallest integer $l \in \N$ such that
every sequence $S$ over $G$ of length $|S| \ge l$ has a nontrivial
zero-sum subsequence;

\item $\mathsf d (G)$ the maximal length of a zero-sum free sequence
over $G$.
\end{itemize}
Then $\mathsf D (G)$ is called the {\it Davenport constant} of $G$,
and we  have trivially that
\[
\mathsf d^* (G) \le \mathsf d (G) = \mathsf D (G) - 1 \,.
\]
We will use without further mention that equality holds for
$p$-groups  and for groups of rank $\mathsf r (G) \le 2$ (\cite[Theorems
5.5.9 and 5.8.3]{Ge-HK06a}) (equality holds for further groups, but
not in general \cite[Corollary 4.2.13]{Ge09a}).

\medskip
{\bf Group algebras and characters.} Let $R$ be a commutative ring \
(throughout, we assume that $R$ has a unit element $1 \ne 0$) \ and
$G$ a finite abelian group. The \ {\it group algebra} \ $R [G]$ \ of
$G$ over $R$ is a free $R$-module with basis $\{ X^g \mid g \in G\}$
\ (built with a symbol $X$), where multiplication is defined by
\[
\Bigl( \sum_{g \in G} a_g X^g \Bigr) \Bigl( \sum_{g \in G} b_g X^g
\Bigr) = \sum_{ g \in G} \Bigl( \sum_{h \in G} a_h b_{g-h} \Bigr)
X^g \,.
\]
We view $R$ as a subset of $R [G]$ by means of $a = a X^0$ for all
$a \in R$. An element of $R$ is a zero-divisor [\,a unit\,] of
$R[G]$ if and only if it is a zero-divisor [\,a unit\,] of $R$.

\smallskip
Let $K$ be a field, \ $G$ a finite abelian group with $\exp(G) = n
\in \mathbb N$ \ and \ $\mu_n(K) = \{\zeta \in K \mid \zeta^n =1\}$ the
group of $n$-th roots of unity in $K$. An $n$-th root of unity $\zeta$
is called \ {\it primitive} \ if $\zeta^m \ne 1$ \ for all $m \in
[1,n-1]$, and we denote by $\mu_n^* (K) \subset \mu_n (K)$  the
subset of all primitive $n$-th roots of unity. We denote by \ ${\rm
Hom} (G, K^{\times}) = {\rm Hom} (G, \mu_n(K))$ \ the \ {\it
character group of \ $G$ \ with values in $K$} (whose operation is given by pointwise multiplication with the constant $1$ function as identity), and we briefly set
$\widehat G = {\rm Hom} (G, K^{\times})$ if there is no danger of
confusion. Every character \ $\chi \in \widehat G$ \ has a unique
extension to a $K$-algebra homomorphism \ $\chi \colon K[G] \to K$
(again denoted by $\chi$) acting by means of
\[
\chi \Bigl( \sum_{g \in G} a_g X^g \Bigr) = \sum_{g \in G} a_g \chi
(g) \,.
\]
We call $K$ a  \ {\it splitting field} \ of $G$ if \ $|\mu_n(K)| =
n$. Let $K$ be a splitting field of $G$ and $\widehat G = {\rm Hom}
(G, K^{\times})$. We gather the properties needed for the sequel (for
details see \cite[Section 5.5]{Ge-HK06a} and \cite[\S 17]{Cu-Re81}). We have \ ${\rm char} (K) \nmid \exp (G)$, \
$|G| = |G| \,1_K \in K^{\times}$,
 \ $G \cong {\rm Hom} (G, K^{\times})$, and
the map
\[
{\rm Hom} (G, K^{\times}) \times G \to K^\times\,, \quad
\text{defined by} \quad (\chi, \, g) \ \mapsto \ \chi(g)\,,
\]
is a non-degenerate pairing (that is, if $\chi (g) = 1$ for all
$\chi \in \widehat G$, then $g= 0$, and if $\chi (g) = 1$ for all $g
\in G$, then clearly $\chi = 1$, the constant $1$ function).

Furthermore,  the Orthogonality Relations hold (\cite[Proposition
5.5.2]{Ge-HK06a}), and for every $f \in K[G]$, we have (see
\cite[Proposition 5.5.2]{Ge-HK06a})
$$f = 0\in K[G]\;\mbox{ if and only if }\chi (f) = 0\;\mbox{ for
every }\;\chi \in {\rm Hom} (G, K^{\times}).$$ Moreover, if $\chi(f)
\ne 0$ for all $\chi \in {\rm Hom}(G,K^\times)$, then $f \in
K[G]^\times$; explicitly, a simple calculation using the
Orthogonality Relations shows that
\[
f^{-1} = \frac 1 {|G|} \sum_{g \in G} \Bigl( \sum_{\chi \in {\rm
Hom}(G,K^\times) } \frac{\chi(-g)}{\chi(f)} \,\Bigr)\, X^g \ \,.
\]
For a subgroup $H \subset G$, we set
\[
H^\perp = \{ \chi \in \widehat G \mid \chi (h) = 1 \ \text{for all}
\ h \in H \} \,.
\]
We clearly have a natural isomorphism $H^\perp \cong \widehat{G/H}$.

\bigskip
\section{Proof of the Theorem}
\bigskip

We fix our notation, which will remain valid throughout this section.
Let $G = \C_m \oplus \C_{mn}$ with $m \in \N_{\ge 2}$, $n \in \N$ and let $e_1,e_2 \in G$ be such that $G = \langle e_1 \rangle \oplus \langle e_2 \rangle$, $\ord(e_1)=m$ and $\ord(e_2)=mn$.
Furthermore, let $K$ be a splitting field of $G$, $\zeta \in \mu^*_{mn}(K)$, and let $\psi,\varphi \in \widehat G$ be defined by $\psi(e_1)=\zeta^n$, $\psi(e_2)=1$ and $\varphi(e_1)=1$, $\varphi(e_2)=\zeta$.
Then $\ord(\psi)=m$, $\ord(\varphi)=mn$ and $\widehat G = \langle\psi\rangle \oplus\langle\varphi\rangle$.

Note that, in the case $m=p \in \PP$,
\[
\theta \colon
\begin{cases}
  \Fp \times \langle\psi,\varphi^n\rangle &\to \langle\psi,\varphi^n\rangle\\
  (k+p\Z, \chi) &\mapsto \chi^k,
\end{cases}
\]
is an $\Fp$-vector space structure on
$(\langle\psi,\varphi^n\rangle, \cdot)$. Whenever
$\langle\psi,\varphi^n\rangle$ is considered as $\Fp$-vector space
it is done so with respect to $\theta$.

The following Lemmas \ref{lemma:reduce-multiples} and \ref{lemma:G0}
will allow us to restrict ourselves to sequences consisting of
certain special elements in the proof of Theorem \ref{thm:main}.1.
Lemma \ref{lemma:G0} is a generalization of a statement used by W.
Gao and Y. Li in their proof of the case $m=2$ \cite{Ga-Li09a}.

\medskip
\begin{lemma} \label{lemma:reduce-multiples}
  Let $R$ be a  commutative ring, $g_1\cdot\ldots\cdot g_l \in \mathcal{F}(G)$ a sequence over $G$, and let $a_1,\ldots,a_l \in R\setminus \{0\}$ be such that $(X^{g_1} - a_1) \cdot\ldots\cdot (X^{g_l} - a_l) = 0 \in R[G]$
  Then, for any $k_1,\ldots,k_l \in \N$,  also $(X^{k_1 g_1} - a_1^{k_1}) \cdot\ldots\cdot (X^{k_l g_l} - a_l^{k_l}) = 0 \in R[G]$.
\end{lemma}

\begin{proof}
  For all $i \in [1,l]$,
  \[
  X^{k_i g_i} - a_i^{k_i} = (X^{g_i} - a_i) \sum_{j=0}^{k_i-1} X^{j g_i} (a_i)^{k_i - 1 -j},
  \]
  from which the lemma immediately follows.
\end{proof}

\medskip
\begin{lemma} \label{lemma:G0}
  Let $R$ be a  commutative ring and
  \[
  G_0 = \{ e_1 \} \cup \Big\{ k e_1 + \prod_{p \in \PP,\,  p \mid m} p^{u_p} e_2 \mid k \in [0,m-1], u_p \in \N_0 \Big\}.
  \]
  Let $M \in \N$ be such that, for every sequence $S=g_1 \cdot\ldots\cdot g_{M+1} \in \mathcal{F}(G_0)$, there exist $a_1,\ldots,a_{M+1} \in R\setminus \{0\}$ such that
  \[
    f = (X^{g_1} - a_1) \cdot\ldots\cdot (X^{g_{M+1}} - a_{M+1}) = 0 \in R[G].
  \]
  Then $\dd(G,R) \le M$.
\end{lemma}

\begin{proof}
  By Lemma \ref{lemma:reduce-multiples} and the definition of $\dd(G,R)$, it is sufficient to show that every element $g \in G$ is a multiple of an element in $G_0$.

  Let $g=k e_1 + l e_2$ with $k \in [0,m-1]$ and $l \in [0,mn-1]$.
  If $l=0$, $g$ is obviously a multiple of $e_1$.
  Consider the case $l \ne 0$.
  Then $l = \prod_{p \in \PP, p \mid m} p^{\val_p(l)} \cdot q$ with $q \in [1,mn-1]$ and $\gcd(q,m)=1$.
  Therefore there exists an $a \in [1,m-1]$ with $qa \equiv 1 \mod m$.
  From $\ord(e_1)=m$, it follows that $g = q(ak e_1 + \prod_{p \in \PP, p \mid m} p^{\val_p(l)} e_2)$.
  Choosing $k' \in [0,m-1]$ such that $k' \equiv ak \mod m$, we obtain $g = q(k' e_1 + \prod_{p \in \PP, p \mid m} p^{\val_p(l)} e_2)$, which is a multiple of an element in $G_0$.
\end{proof}

\medskip
\begin{lemma} \label{lemma:orthogonal}
  Let $g \in G$ and $\chi, \chi' \in \widehat G$. Then $\chi'(g)=\chi(g)$ if and only if $\chi' \in \chi \langle g\rangle^\perp$. Also
  \begin{enumerate}
    \item \label{el:ke1+e2} $\langle k e_1 + e_2 \rangle^\perp = \langle \psi\varphi^{-nk} \rangle$ for $k \in [0,m-1]$;
    \item \label{el:ke1+mle2} $\langle \varphi^n \rangle\subset\langle k e_1 + m l e_2 \rangle^\perp$ for $k \in [0,m-1]$ and $l \in [0,n-1]$.
  \end{enumerate}
\end{lemma}

\begin{proof}
  Clearly $\chi'(g) = \chi(g)$ if and only if $\chi^{-1} \chi'(g) = 1$, i.e., $\chi' \in \chi \langle g \rangle^\perp$.

\smallskip
1. From $\psi^{-1}(k e_1 + e_2) = \zeta^{-nk} = \varphi^{-nk}(k e_1
+ e_2)$, it follows that $\langle \psi\varphi^{-nk} \rangle \subset
\langle k e_1 + e_2 \rangle^\perp$.
      Then $\ord(k e_1 + e_2) = mn$ and $\langle k e_1 + e_2 \rangle^\perp \cong \widehat{ G/\langle k e_1+e_2 \rangle }$ imply $\card{\langle k e_1 + e_2 \rangle^\perp}=m$, from which $\langle k e_1 + e_2 \rangle^\perp = \langle \psi\varphi^{-nk} \rangle$ follows.

\smallskip
2. Observe that $\varphi^n(ke_1 + ml e_2) = \zeta^{nml} = (\zeta^{nm})^l = 1$
implies $\langle \varphi^n \rangle \subset \langle k e_1 + ml e_2
\rangle^\perp$.
%
\end{proof}

\medskip
\begin{lemma} \label{lemma:cover-equiv}
  Let $H \subset \widehat G$ and $S=g_1\cdot\ldots\cdot g_l \in \mathcal{F}(G)$.
  Then the following statements are equivalent\,{\rm :}
  \begin{enumerate}
    \item[(a)] There exist $a_1,\ldots,a_l \in K^\times$ such that $\chi\big( \prod_{i=1}^l (X^{g_i}-a_i)\big)=0$ for all $\chi \in H$.
    \item[(b)] There exist $s \in [0, l]$ and $\chi_1,\ldots,\chi_s \in H$ such that  $H \subset \bigcup_{i=1}^s \chi_i \langle g_i \rangle^\perp$.
    \item[(c)] $H=\emptyset$ or there exist $\chi_1,\ldots,\chi_l \in H$ such that $H \subset \bigcup_{i=1}^l \chi_i \langle g_i \rangle^\perp$.
  \end{enumerate}
\end{lemma}

\begin{proof}
  For $H=\emptyset$ all statements are trivially true.
  Let $H \neq \emptyset$.

  (a) $\Rightarrow$ (b)
  The extension of $\chi \in \widehat G$ to $K[G]$ is a $K$-algebra homomorphism, and thus $$\chi\big( \prod_{i=1}^l (X^{g_i}-a_i) \big) =0$$ if and only if there is an $i \in [1,l]$ with $\chi(X^{g_i}-a_i)=0$, i.e., $\chi(g_i)=a_i$.
  Let
  \[
  s = \card{ \{ i \in [1,l] \mid \text{ there exists a $\chi \in H$ such that } \chi(g_i) = a_i \} } \in [0,l].
  \]
  Without restriction let $g_1,\ldots,g_s$ and $a_1,\ldots,a_s$ be such that there exist $\chi_i \in H$ with $\chi_i(g_i)=a_i$ for $i \in [1,s]$.
  Let $\chi \in H$.
  Then, by assumption, $\chi(g_i)=a_i$ for some $i \in [1,s]$.
  Therefore $\chi_i^{-1} \chi(g_i) = 1$, i.e. $\chi \in \chi_i \langle g_i\rangle^\perp$.

  (b) $\Rightarrow$ (a)
  Let $a_i=\chi_i(g_i)$ for $i \in [1,s]$ and let $a_{s+1}=\ldots=a_l=1$.
  Let $\chi \in H$.
  Then, by assumption, there exists an $i \in [1,s]$ such that $\chi \in \chi_i \langle g_i \rangle^\perp$, i.e.,
  $\chi(g_i)=\chi_i(g_i)=a_i$. Hence $\chi(X^{g_i}-a_i)=0$.

  (b) $\Leftrightarrow$ (c) Obvious.
\end{proof}

\medskip
Note that, in particular, $\dd(G,K)$ is the supremum of all $l \in \N_0$ such that there exists a sequence $S=g_1\cdot\ldots\cdot g_l \in \mathcal{F}(G)$ with
\[
    \bigcup_{i=1}^l \chi_i \langle g_i \rangle^\perp \subsetneq \widehat G
\]
for any choice of $\chi_1,\ldots,\chi_l \in \widehat G$.
Or, equivalently, $\dd(G,K)+1$ is the minimum of all $l \in \N_0$ such that, for any sequence $S=g_1\cdot\ldots\cdot g_l \in \mathcal{F}(G)$, there exist $\chi_1,\ldots,\chi_l \in \widehat G$ such that $\widehat G$ can be covered as above:
\[
    \widehat G = \bigcup_{i=1}^l \chi_i \langle g_i \rangle^\perp.
\]
Consider $m=p \in \PP$.
Our strategy for finding an upper bound on $\dd(G,K)$ will be to subdivide $\widehat G$ into cosets modulo $\langle\psi,\varphi^n\rangle$ and cover each of these cosets individually.
Lemma \ref{lemma:G0} allows us to restrict ourselves to certain special elements $g \in G$ in doing so, and from Lemma \ref{lemma:orthogonal}, we see that for these elements $\langle g \rangle^\perp$ contain (or are) $1$-dimensional subspaces, i.e., lines of the $2$-dimensional $\Fp$-vector space $\langle \psi,\varphi^n \rangle$. Then, for $\chi \in \langle\psi,\varphi^n\rangle$, $\chi \langle g\rangle^\perp$ is an affine line in $\langle \psi,\varphi^n \rangle$ containing the ``point'' $\chi$, and our task essentially boils down to covering $n$ copies of $\langle\psi,\varphi^n\rangle$ by such lines (where the slopes are fixed by $S$).

Before we do so, we study some simple configurations in Lemma \ref{lemma:parallel} and Lemma \ref{lemma:star}.
The main part of the proof for the cases $m \in \{2,3\}$ then follows in Lemma \ref{lemma:residue-seqs}.
It is based on the proof by Gao and Li of the case $m=2$, but is stated in terms of group characters instead of working with the group algebra directly.

\medskip
\begin{lemma} \label{lemma:parallel}
  Let $s \in [0,m]$ and let $S=g_1\cdot\ldots\cdot g_{s+(m-s)m} \in \mathcal{F}(G)$ such that either $g_1=\ldots=g_s=ke_1+e_2$ with $k \in [0,m-1]$ or $g_1,\ldots,g_s \in \{ k e_1 + m l e_2 \mid k \in [0,m-1], l \in \N_0 \}$.
  Then there exist $\chi_1,\ldots,\chi_{s+(m-s)m}$ such that $\langle\psi,\varphi^n\rangle \subset \bigcup_{i=1}^{s+(m-s)m} \chi_i \langle g_i \rangle^\perp$.
\end{lemma}

\begin{proof}
  Let $L=\langle\psi\varphi^{-nk}\rangle$ in the case $g_1=\ldots=g_s=ke_1+e_2$, and let $L=\langle\varphi^n\rangle$ otherwise.
  Since $L$ is a subgroup of $\langle\psi,\varphi^n\rangle$ and has cardinality $\card{L}=m$, there exist $\tau_1,\ldots,\tau_m \in \langle\psi,\varphi^n\rangle$ such that $\langle\psi,\varphi^n\rangle = \biguplus_{i=1}^m \tau_i L$.
  By Lemma \ref{lemma:orthogonal}, $L \subset \langle g_i \rangle^\perp$ for $i \in [1,s]$.
  Then
  \[
  \langle\psi,\varphi^n\rangle \; \subset \; \bigcup_{i=1}^s \tau_i \langle g_i \rangle^\perp \cup \biguplus_{i=s+1}^m \tau_i L.
  \]

  For $j \in [s+1,s+(m-s)m]$, let $\chi'_j \in \langle g_j \rangle^\perp$, and let $L=\{\lambda_1,\ldots,\lambda_{m}\}$.
  Then, for $i \in [s+1,m]$,
  \[
  \tau_i L = \{ \tau_i \lambda_j \mid j \in [1,m] \} \subset \bigcup_{j=1}^{m} \tau_i \lambda_j \chi_{s+(i-(s+1))m+j}'^{-1} \langle g_{s+(i-(s+1))m+j}
  \rangle^\perp \,. \qedhere
  \]
\end{proof}

\medskip
\begin{lemma} \label{lemma:star}
  Let $m=p \in \PP$, $g \in \{ k e_1 + p l e_2 \mid k \in [0,p-1], l \in \N_0 \}$ and $S=\prod_{i=0}^{p-1} (i e_1 + e_2) g$.
  Then $\langle\psi,\varphi^n\rangle \subset \bigcup_{i=0}^{p-1} \langle i e_1 + e_2\rangle^\perp \cup \langle g \rangle^\perp$.
\end{lemma}

\begin{proof}
  By Lemma \ref{lemma:orthogonal},
  \[
  \bigcup_{i=0}^{p-1} \langle\psi\varphi^{-ni}\rangle \cup \langle\varphi^n\rangle \; \subset \; \bigcup_{i=0}^{p-1} \langle i e_1 + e_2\rangle^\perp \cup \langle g \rangle^\perp.
  \]
  Let $\psi^k\varphi^{nl} \in \langle\psi,\varphi^n\rangle$ with $k,l\in[0,p-1]$.
  In the case $k=0$, clearly $\varphi^{nl} \in \langle\varphi^n\rangle$.
  Otherwise, there exists an $i \in [0,p-1]$ such that $-ik \equiv l \mod p$.
  Hence $\psi^k\varphi^{nl} = (\psi\varphi^{-ni})^k \in \langle\psi\varphi^{-ni}\rangle$.
\end{proof}

\medskip
\begin{lemma} \label{lemma:residue-seqs}
  Let $m=p \in \PP$, $G_1 = \{ e_1 \} \cup \{ k e_1 + p^u e_2 \mid k \in [0,p-1], u \in \N \}$, and
  \[
  G_0 = \{ e_1 \} \cup \{ k e_1 + p^u e_2 \mid k \in [0,p-1], u \in \N_0 \} = \{ k e_1 + e_2 \mid k \in [0,p-1] \} \uplus G_1.
  \]
  If, for all sequences $T=h_1 \cdot\ldots\cdot h_{rp-1} \in \mathcal{F}(G_0)$ with $r \in [2,\min{ \{ p-1, n+1 \} }]$ and $\val_g(T) < p$ for all $g \in G_0$ as well as $\sum_{g \in G_1} \val_g(T) < p$, there exist $\chi_1,\ldots,\chi_{rp-1} \in \widehat G$ such that $\bigcup_{i=0}^{r-2} \varphi^i \langle\psi,\varphi^n\rangle \subset \bigcup_{i=1}^{rp-1} \chi_i \langle h_i \rangle^\perp$, then $\dd(G,K)=\dd^*(G)$.
\end{lemma}

\begin{proof}
  Since $\dd^*(G) \le \dd(G) \le \dd(G,K)$ always holds, it is sufficient to show that $\dd(G,K) \le \dd^*(G) = (pn-1) + (p-1) = (n+1)p - 2$.
  By Lemma \ref{lemma:G0}, it is sufficient to show that, for any sequence $S=g_1\cdot\ldots\cdot g_{(n+1)p -1} \in \mathcal F(G_0)$, there exist $a_1,\ldots,a_{(n+1)p-1} \in K^\times$ such that
  \[
    f = \prod_{i=1}^{(n+1)p-1} (X^{g_i} - a_i) = 0 \in K[G].
  \]
  To see this, we use Lemma \ref{lemma:cover-equiv} and show that there exist $\chi_1,\ldots, \chi_{(n+1)p-1}$ such that
  \[
  \widehat G = \biguplus_{i=0}^{n-1} \varphi^i \langle\psi,\varphi^n\rangle \subset \bigcup_{i=1}^{(n+1)p-1} \chi_i \langle g_i \rangle^\perp.
  \]

  We group the elements of $S$ into as many $p$-tuples of the forms $(e_2,\ldots,e_2)$, $(e_1 + e_2,\ldots, e_1 + e_2)$, \ldots, $( (p-1)e_1 + e_2, \ldots (p-1) e_1 + e_2)$ and $(g'_1,\ldots,g'_p) \in G_1^p$ as possible to obtain $l \in [0,n]$ such tuples.
  Without restriction, let these $p$-tuples be $(g_1,\ldots,g_p)$, \ldots, $(g_{(l-1)p+1},\ldots, g_{lp})$.

  For each $i \in [1,l]$, the tuple $(g_{(i-1)p+1},\ldots,g_{ip})$ fulfills the conditions of Lemma \ref{lemma:parallel} with $s=p$.
  Therefore, there exist $\chi_{(i-1)p+1},\ldots,\chi_{ip}$ such that
  \[
  \varphi^{n-i} \langle\psi,\varphi^n\rangle \; \subset\;  \bigcup_{j=(i-1)p+1}^{ip} \chi_{j} \langle g_{j} \rangle^\perp.
  \]

  It remains to be shown that $\chi_{lp+1},\ldots,\chi_{(n+1)p-1}$ can be chosen such that
  \[
  \bigcup_{i=0}^{n-l-1} \varphi^i \langle\psi,\varphi^n\rangle\; \subset\; \bigcup_{j=lp+1}^{(n+1)p-1} \chi_{j}\langle g_{j}\rangle^\perp.
  \]
  In the case $l \geq n$, this is trivially so, and therefore it is sufficient to consider $l \leq n-1$.

  By $T=g_{lp+1}\cdot\ldots\cdot g_{(n+1)p-1}$ we denote the subsequence of $S$ consisting of the remaining elements.
  We have $\length{T} = \length{S} - lp = (n+1-l)p - 1$.
  In the process of creating $p$-tuples, we partitioned the elements of $G_0$ into $p+1$ different types.
  If there were at least $p$ elements of one type, we could create another tuple, in contradiction to the maximal choice of $l$.
  Thus we must have $\val_g(T) < p$ for all $g \in G_0$, $\sum_{g \in G_1} \val_g(T) < p$, and $\length{T} \le (p+1)(p-1) = p^2 - 1$, which implies $n+1-l \le p$.

  Altogether, we have $n+1-l \in [2,p]$.
  In the case $n+1-l \le p-1$, we set $r=n+1-l \in [2,\min{\{p-1,n+1\}}]$.
  Then, by assumption, $\chi_{lp+1},\ldots,\chi_{(n+1)p-1}$ can be chosen such that $$\bigcup_{i=0}^{r-2} \varphi^i \langle \psi, \varphi^n \rangle\; \subset\; \bigcup_{j=lp+1}^{(n+1)p-1} \chi_{j}\langle g_{j}\rangle^\perp.$$
  Since $r-2 = n-l-1$, this already means $\widehat G \subset \bigcup_{i=1}^{(n+1)p-1} \chi_i \langle g_i \rangle^\perp$.

  In the case $n+1-l = p$, we have $\card{T}=p^2 - 1=(p+1)(p-1)$.
  This can only happen if each of the $p+1$ different types of elements occurs exactly $p-1$ times.
  Therefore $$T=\prod_{j=0}^{p-1} (je_1+e_2)^{p-1} \cdot \prod_{i=0}^{p-2} h_j = \prod_{i=0}^{p-2} \Big( \prod_{j=0}^{p-1} (j e_1 + e_2) \cdot h_i \Big)$$ with $h_0,\ldots,h_{p-2} \in G_1$.
  Without restriction, for $i \in [0,p-2]$, let $g_{lp+i(p+1)+1}\cdot\ldots\cdot g_{lp+i(p+1)+(p+1)} = \prod_{j=0}^{p-1} (j e_1 + e_2) \cdot h_i$.
  For every $i \in [0,p-2]$, we set $\chi_{lp+i(p+1)+1}=\ldots=\chi_{lp+i(p+1)+(p+1)}=\varphi^i$.
  Then, from Lemma \ref{lemma:star}, it follows that $\varphi^i \langle \psi,\varphi^n \rangle \subset \bigcup_{j=i(p+1)+1}^{i(p+1)+(p+1)} \chi_{lp+j} \langle g_{lp+j} \rangle^\perp$.
  Due to $n-l-1=p-2$, this again implies $\widehat G \subset \bigcup_{i=1}^{(n+1)p-1} \chi_i \langle g_i \rangle^\perp$.
\end{proof}

\bigskip
\begin{proof}[\textbf{\emph{Proof of Theorem \ref{thm:main}.1}}]
  For $p=2$, i.e. $G = \C_2 \oplus \C_{2n}$, this follows trivially from Lemma \ref{lemma:residue-seqs}, since there are no admissible sequences.

  Consider $p=3$, i.e., $G = \C_3 \oplus \C_{3n}$.
  Let $G_1 = \{ e_1 \} \cup \{ k e_1 + 3^u e_2 \mid k \in [0,2], u \in \N \}$ and $G_0 = \{ e_2, e_1+e_2, 2e_1 + e_2 \} \uplus G_1$.
  Then, by Lemma \ref{lemma:residue-seqs}, it is sufficient to show that, for $T=h_1\cdot\ldots \cdot h_5 \in \mathcal{F}(G_0)$, we can choose $\chi_1,\ldots,\chi_5 \in \widehat G$ such that $\langle \psi,\varphi^n \rangle \; \subset \; \chi_1 \langle h_1 \rangle^\perp \cup \ldots \cup \chi_5 \langle h_5 \rangle^\perp$.
  We divide the elements into four types: $e_2$, $e_1+e_2$, $2e_1+e_2$ and elements from $G_1$.
  Since $\length{T}=5$, one of these types must occur at least twice.
  Without restriction, let $h_1$ and $h_2$ be of the same type.
  Thus we have either $h_1=h_2=k e_1+e_2$ for some $k \in [0,2]$ or $h_1,h_2 \in G_1$.
  Then $T$ fulfills the conditions of Lemma \ref{lemma:parallel} with $s=2$, and it follows that $\chi_1,\ldots,\chi_5$ can be chosen such that $\langle \psi,\varphi^n\rangle \subset \bigcup_{i=1}^5 \chi_i \langle h_i \rangle^\perp$.
\end{proof}

The following Lemma \ref{lemma:lines} recapitulates a few simple facts, which are well known in the context of affine lines, and will be used extensively in the construction of a counterexample in the case $p \ge 5$ and $n \ge 2$.

\medskip
\begin{lemma} \label{lemma:lines}
  Let $m=p \in \PP$,  $g_1 = k_1 e_1 + e_2$, $g_2 = k_2 e_1 + e_2$ with $k_1,k_2 \in
  [0,p-1]$,
   $\chi \in \widehat G$ and $\chi_1,\chi_2 \in \chi \langle\psi,\varphi^n\rangle$.

  \begin{enumerate}
    \item $\chi^{-1} \chi_i \langle g_i \rangle^\perp = \varphi^{n s_i} \langle g_i \rangle^\perp$ with $s_i \in [0,p-1]$ for $i \in \{1,2\}$.

    \item $\chi^{-1} \chi_i \langle g_i \rangle^\perp = \{ \psi^u \varphi^{nv} \mid u,v \in [0,p-1] \text{ with } k_i u + v \equiv s_i \mod p \}$ for $i \in \{1,2\}$.

    \item
      \begin{enumerate}
        \item $\card{\chi_1 \langle g_1 \rangle^\perp \cap \chi_2 \langle g_2 \rangle^\perp} = 1$ if and only if $g_1 \neq g_2$.
        \item $\card{\chi_1 \langle g_1 \rangle^\perp \cap \chi_2 \langle g_2 \rangle^\perp} = 0$ if and only if $g_1 = g_2$ and $s_1 \neq s_2$.
        \item $\card{\chi_1 \langle g_1 \rangle^\perp \cap \chi_2 \langle g_2 \rangle^\perp} = p$ if and only if $g_1 = g_2$ and $s_1 = s_2$.
      \end{enumerate}
  \end{enumerate}
\end{lemma}

\begin{proof}
1. Let $i \in \{1,2\}$ and $\chi^{-1} \chi_i=\psi^{u_i}\varphi^{n
v_i}$ with $u_i,v_i \in [0,p-1]$.
      By Lemma \ref{lemma:orthogonal}.1, $\langle g_i \rangle^\perp = \langle\psi\varphi^{-n k_i} \rangle$.
      Therefore $\varphi^{-n(k_i u_i + v_i)} \chi^{-1} \chi_i = \psi^{u_i} \varphi^{-n k_i u_i} \in \langle g_i \rangle^\perp$, and hence $\chi^{-1} \chi_i \langle g_i \rangle^\perp = \varphi^{n s_i} \langle g_i \rangle^\perp$ with $s_i \in [0,p-1]$ chosen such that $s_i \equiv k_i u_i + v_i \mod p$.

\smallskip
2. In view of Lemma \ref{lemma:orthogonal}.1, we have, for $u,v \in [0,p-1]$, \, $\psi^u \varphi^{nv} \in \chi^{-1}
\chi_i \langle g_i \rangle^\perp = \varphi^{n s_i} \langle \psi
\varphi^{-n k_i}\rangle$ if and only if $\psi^u \varphi^{nv} =
\psi^w \varphi^{n(s_i - k_i w)}$ for some $w \in [0,p-1]$.
      This is the case if and only if $u \equiv w \mod p$ and $v \equiv s_i - k_i w \mod p$, i.e., if and only if $u \equiv w \mod p$ and $k_i u + v \equiv s_i \mod p$ (recall by Lemma \ref{lemma:orthogonal}.1 that $\langle g_i\rangle^ \perp\subset \langle \psi,\varphi^n\rangle$).

\smallskip
3. By 2, we have $\chi^{-1} \chi_1 \langle g_1 \rangle^\perp \cap
\chi^{-1} \chi_2 \langle g_2 \rangle^\perp = \{ \psi^u \varphi^{nv}
\mid u,v \in [0,p-1] \text{ with } k_1 u + v \equiv s_1 \mod p
\text{ and } k_2 u + v \equiv s_2 \mod p \}$.
      Since
      \[
        \card{\chi^{-1} \chi_1 \langle g_1 \rangle^\perp \cap \chi^{-1} \chi_2 \langle g_2 \rangle^\perp} = \card{\chi_1 \langle g_1 \rangle^\perp \cap \chi_2 \langle g_2 \rangle^\perp},
      \]
      it is sufficient to consider the number of solutions of the linear system
      \[
        k_1 u + v \equiv s_1 \mod p \qquad\textrm{ and }\qquad k_2 u + v \equiv s_2 \mod p
      \]
      for $u,v\in[0,p-1]$ over $\Fp$.
      In the case $g_1 \neq g_2$, i.e., $k_1 \ne k_2$, it possesses a unique solution.
      In the case $g_1 = g_2$, it possesses no solution for $s_1 \neq s_2$.
      For $s_1=s_2$, the two equations coincide, and we obtain $p$ solutions.
\end{proof}

In the construction of the counterexamples, we use the same characterization of $\dd(G,K)$, derived from Lemma \ref{lemma:cover-equiv}, as in the proof of Theorem \ref{thm:main}.1---except now we show that it is not possible to cover $\widehat G$ with such subsets.
To do so, we first consider a special type of sequence in Lemma \ref{lemma:l-triple}, which will turn out to be the only one which cannot be discarded with simpler combinatorial arguments, as will be given in the proof of Theorem \ref{thm:main}.2 that follows the lemma.

\medskip
\begin{lemma} \label{lemma:l-triple}
  Let $m=p \in \PP$, $p \ge 5$ and $k_1,k_2,k_3 \in [0,p-1]$ be distinct.
  Let $l \in [2,p-1]$,
  \[
  T = (k_1 e_1 + e_2)^l (k_2 e_1 + e_2)^l (k_3 e_1 + e_2)^l \in \mathcal{F}(G),
  \]
  and $\chi \in \widehat G$.
  For $i \in [1,3]$ and $j \in [1,l]$, let $\chi_{i,j} \in \widehat G$.
  Then
  \[
  \lrcard{ \Big( \bigcup_{i=1}^3 \bigcup_{j=1}^l \chi_{i,j} \langle k_i e_1+e_2 \rangle^\perp \Big) \cap \chi \langle \psi,\varphi^n\rangle } < l (3p - 2l).
  \]
\end{lemma}

\begin{proof}
  We set $g_i=k_i e_1 + e_2$ for $i \in [1,3]$.
  Let $i \in [1,3]$ and $j \in [1,l]$.
  We can assume $\chi_{i,j} \in \chi \langle\psi,\varphi^n\rangle$ since otherwise $\chi_{i,j} \langle g_i\rangle^\perp \cap \chi \langle\psi,\varphi^n\rangle = \emptyset$ (due to $\langle g_i\rangle^\perp = \langle\psi\varphi^{-nk_i} \rangle \subset \langle\psi,\varphi^n\rangle$).
  Using Lemma \ref{lemma:lines}.1, we can furthermore assume $\chi^{-1} \chi_{i,j}=\varphi^{ns_{i,j}}$ with $s_{i,j} \in [0,p-1]$.
  And we can then also assume, without restriction, $s_{i,j} \neq s_{i,j'}$ for $j' \in [1,l]\setminus \{j\}$, since otherwise $\chi_{i,j}\langle g_i\rangle^\perp = \chi_{i,j'} \langle g_i \rangle^\perp$.

  For $i \in [1,3]$, let $E_i = \bigcup_{j=1}^l \chi_{i,j} \langle g_i \rangle^\perp$.
  Then
  \[
  \Big( \bigcup_{i=1}^3 \bigcup_{j=1}^l \chi_{i,j} \langle g_i \rangle^\perp \Big) \cap \chi \langle \psi,\varphi^n\rangle = E_1 \cup E_2 \cup E_3
  \]
  and
  \[
  \card{E_1 \cup E_2 \cup E_3} = \sum_{i=1}^3 \card{E_i} - \sum_{1 \le i < i' \le 3} \card{E_i \cap E_{i'}} + \card{E_1 \cap E_2 \cap E_3}.
  \]
  For $i,i' \in [1,3]$ distinct, we show $\card{E_i}=lp$, $\card{E_i \cap E_{i'}} = l^2$ and $\card{E_1 \cap E_2 \cap E_3} < l^2$.
  Then $\card{E_1 \cup E_2 \cup E_3} < 3lp - 3l^2 + l^2 = l(3p - 2l)$.

  Let $i \in [1,3]$.
  By Lemma \ref{lemma:lines}.3b, $\chi_{i,j}\langle g_i \rangle^\perp \cap \chi_{i,j'}\langle g_i\rangle^\perp = \emptyset$ for $j,j' \in [1,l]$ with $j \neq j'$, and $\card{\langle g_i \rangle^\perp}=\card{\langle \psi\varphi^{-n k_i}\rangle} = p$ (by Lemma \ref{lemma:orthogonal}.1).
  Therefore $\card{E_i}=lp$.

  Let $i,i' \in [1,3]$ be distinct.
  For $j,j' \in [1,l]$ distinct, we have $\chi_{i,j} \langle g_i \rangle^\perp \cap \chi_{i,j'} \langle g_i \rangle^\perp = \emptyset$ and $\chi_{i',j} \langle g_{i'} \rangle^\perp \cap \chi_{i',j'} \langle g_{i'} \rangle^\perp = \emptyset$ (by Lemma \ref{lemma:lines}.3b).
  This implies that, for
  \[
  E_i \cap E_{i'} = \big( \bigcup_{j=1}^l \chi_{i,j} \langle g_i\rangle^\perp \big) \cap
                    \big( \bigcup_{j'=1}^l \chi_{i',j'} \langle g_{i'}\rangle^\perp \big)
                  = \biguplus_{j=1}^l \biguplus_{j'=1}^l ( \chi_{i,j} \langle g_i\rangle^\perp \cap \chi_{i',j'} \langle g_{i'} \rangle^\perp ),
  \]
  the union is disjoint.
  By Lemma \ref{lemma:lines}.3a $\card{\chi_{i,j} \langle g_i\rangle^\perp \cap \chi_{i',j'} \langle g_{i'} \rangle^\perp} = 1$ for $j,j' \in [1,l]$, and therefore $\card{E_i \cap E_{i'}} = l^2$.

  Assume $\card{E_1 \cap E_2 \cap E_2} \ge l^2$.
  Then, since $\card{E_1 \cap E_2} = l^2$, $\card{E_1 \cap E_2 \cap E_3} = l^2$.
  For $a \in \Z$, let $\lbar a = a + p\Z \in \Fp$.
  Let $u,v \in [0,p-1]$.
  By Lemma \ref{lemma:lines}.2, $\chi \psi^u\varphi^{nv} \in E_1 \cap E_2 \cap E_3$ if and only if there are $b_i \in \{ s_{i,1}, \ldots, s_{i,l} \}$, for $i \in [1,3]$, such that
  \begin{align*}
    \lbar{k_1} \lbar u + \lbar v &= \lbar{b_1} \\
    \lbar{k_2} \lbar u + \lbar v &= \lbar{b_2} \\
    \lbar{k_3} \lbar u + \lbar v &= \lbar{b_3}.
  \end{align*}

  Since $\lbar{k_1}, \lbar{k_2}$ and $\lbar{k_3}$ are pairwise distinct, $(\lbar{k_1}, \lbar 1)$, $(\lbar{k_2},\lbar 1)$ and $(\lbar{k_3}, \lbar 1)$ are pairwise $\Fp$-linearly independent.
  For $i \in [1,3]$, we define $\Phi_i: \chi \langle\psi,\varphi^n\rangle \to \Fp$ by $\Phi_i(\chi \psi^u\varphi^{nv}) = \lbar{k_i} \lbar{u} + \lbar{v}$.
  Then the linear independence of $(\lbar{k_1}, \lbar 1)$ and $(\lbar{k_2},\lbar 1)$ implies that $\Phi=(\Phi_1,\Phi_2): \chi \langle\psi,\varphi^n\rangle \to \Fp^2$ is bijective.
  We have $\Phi(E_1\cap E_2\cap E_3) \subset \{ \lbar{s_{1,1}},\ldots,\lbar{s_{1,l}} \} \times \{ \lbar{s_{2,1}}, \ldots, \lbar{s_{2,l}} \}$, and due to $l^2 = \card{E_1 \cap E_2 \cap E_3} \le \card{\{ \lbar{s_{1,1}},\ldots,\lbar{s_{1,l}} \} \times \{ \lbar{s_{2,1}}, \ldots, \lbar{s_{2,l}} \} } = l^2$, equality holds.
  In particular, $\Phi_1(E_1\cap E_2\cap E_3)=\{ \lbar{s_{1,1}},\ldots,\lbar{s_{1,l}} \}$ and $\Phi_2(E_1\cap E_2\cap E_3)=\{ \lbar{s_{2,1}},\ldots,\lbar{s_{2,l}} \}$.

  Because $(\lbar{k_1}, \lbar 1)$, $(\lbar{k_2},\lbar 1)$ and $(\lbar{k_3}, \lbar 1)$ are pairwise $\Fp$-linearly independent, there exist $x,y \in \Fp^\times$ such that $(\lbar{k_3}, \lbar 1) = x (\lbar{k_1}, \lbar 1) + y (\lbar{k_2},\lbar 1)$.
  Hence $\Phi_3 = x \Phi_1 + y \Phi_2$.
  Now
  $\card{x \Phi_1(E_1\cap E_2\cap E_3)} = \card{y \Phi_2(E_1\cap E_2\cap E_3)} = l$. Also, since $x,\,y\neq 0$, we have (similar to $\Phi$) that $(x\Phi_1,y\Phi_2): \chi \langle\psi,\varphi^n\rangle \to \Fp^2$ is a bijective map.  Thus, in view of $\card{x \Phi_1(E_1\cap E_2\cap E_3)} = \card{y \Phi_2(E_1\cap E_2\cap E_3)} = l$ and $|E_1\cap E_2\cap E_3|=l^2$, we see that $(x\Phi_1, y\Phi_2)(E_1 \cap E_2 \cap E_3) = x \Phi_1(E_1\cap E_2\cap E_3) \times y \Phi_2(E_1 \cap E_2 \cap E_3)$.
  Therefore
  \[
  \Phi_3(E_1\cap E_2\cap E_3)=x\Phi_1(E_1\cap E_2\cap E_3)+y\Phi_2(E_1\cap E_2\cap E_3),
  \] where the inclusion ``$\subset$'' is obvious and ``$\supset$'' follows since for any $\alpha,\beta \in E_1 \cap E_2 \cap E_3$ we can find $\theta \in E_1 \cap E_2 \cap E_3$ such that $(x\Phi_1(\alpha), y\Phi_2(\beta)) = (x\Phi_1(\theta), y\Phi_2(\theta))$, and hence in particular $x\Phi_1(\alpha) + y\Phi_2(\beta) = x\Phi_1(\theta) + y\Phi_2(\theta) = \Phi_3(\theta)$.
  From the Cauchy-Davenport Theorem (Lemma \ref{2.1}), it then follows that $\card{\Phi_3(E_1\cap E_2\cap E_3) } \ge \min{ \{ 2l-1, p \} } > l$, a contradiction, since $\Phi_3(E_1\cap E_2\cap E_3) \subset \{ \lbar{s_{3,1}},\ldots,\lbar{s_{3,l}} \}$.
\end{proof}

\bigskip
\begin{proof}[\textbf{\emph{Proof of Theorem \ref{thm:main}.2}}]
  Consider $m=p \in \PP_{\geq 5}$ and $n \ge 2$.
  Let $k_1,\ldots,k_4 \in [0,p-1]$ be pairwise distinct and set $g_i = k_i e_1 + e_2 \in G$ for $i \in [1,4]$.
  Furthermore, set $m_1 = (n-2)p+(p-1)$, $m_2=m_3=p-1$ and $m_4=2$.
  We consider the sequence
  \[
    S = \prod_{i=1}^4 g_i^{m_i} \in \mathcal{F}(G)
  \] and, for any choice of
 $\chi_{i,j} \in \widehat G$ for $i \in [1,4]$ and $j \in [1,m_i]$, show that
  \[
    \bigcup_{i=1}^4 \bigcup_{j=1}^{m_i} \chi_{i,j} \langle g_i \rangle^\perp \subsetneq \widehat G.
  \]
  Then, by Lemma \ref{lemma:cover-equiv} and the definition of $\dd(G,K)$,
  \[
    \dd(G,K) \geq \length{S} = p + pn -1 > p + pn - 2 = \dd^*(G).
  \]

  Let $\chi_{i,j} \in \widehat G$ for $i \in [1,4]$ and $j \in [1,m_i]$ be arbitrary. Assume, to the contrary, $\bigcup_{i=1}^4 \bigcup_{j=1}^{m_i} \chi_{i,j} \langle g_i \rangle^\perp = \widehat G$.
  For $i \in [1,4]$ and $j,j' \in [1,m_i]$ distinct, we can without restriction assume $\chi_{i,j} \langle g_i \rangle^\perp \ne \chi_{i,j'} \langle g_i \rangle^\perp$.

  For any permutation $\sigma \in \mathfrak{S}_n$ (which will be fixed later),
  \[
    \widehat G = \biguplus_{\nu=1}^n \varphi^{\sigma(\nu)} \langle\psi,\varphi^n\rangle.
  \]
  For given $i \in [1,4]$ and $j \in [1,m_i]$, we have by Lemma \ref{lemma:orthogonal} that $\chi_{i,j} \langle g_i \rangle^\perp \subset \varphi^{\sigma(\nu)} \langle\psi,\varphi^n\rangle$ for a uniquely determined $\nu \in [1,n]$.
  For $i \in [1,4]$ and $\nu \in [1,n]$, we can therefore define
  \[
  B_i^{(\nu)} = \big\{ \chi_{i,j} \mid j \in [1,m_i] \text{ with } \chi_{i,j} \langle g_i \rangle^\perp \subset \varphi^{\sigma(\nu)} \langle\psi,\varphi^n\rangle \big\}.
  \]
  We also define $n^{(\nu)} = \max{ \{ \card{B_i^{(\nu)}} \mid i \in [1,4] \} }$ as well as $l^{(\nu)} = \sum_{i=1}^4 \card{ B_i^{(\nu)} }$, for $\nu \in [1,n]$.

  Let $\nu \in [1,n]$.
  By assumption,
  \[
    \varphi^{\sigma(\nu)} \langle\psi,\varphi^n\rangle = \bigcup_{i=1}^4 \bigcup_{\chi \in B_i^{(\nu)}} \chi \langle g_i \rangle^\perp.
  \]
  Thus, since $\card{ \langle\psi,\varphi^n\rangle } = p^2$ and $\card{\langle g_i \rangle^\perp}=p$ for all $i \in [1,4]$, we have $l^{(\nu)} \ge p$.
  On the other hand, $n^{(\nu)} \le p$ because otherwise there would exist $i \in [1,4]$ and $j, j' \in [1,m_i]$ distinct such that $\chi_{i,j} \langle g_i \rangle^\perp \cap \chi_{i,j'} \langle g_i \rangle^\perp \neq \emptyset$,
  but this would already imply $\chi_{i,j} \langle g_i \rangle^\perp = \chi_{i,j'} \langle g_i \rangle^\perp$, contrary to assumption.

 Fix $\sigma \in \mathfrak{S}_n$ so that there is a $k \in \N_0$   such that $n^{(1)},\ldots,n^{(k)} < p$ and $n^{(k+1)}=\ldots=n^{(n)}=p$.
 Since $m_i<p$ for $i\geq 2$, we see (for $\nu \in [1,n]$) that $n^{(\nu)} = p$ is only possible if $\card{B_1^{(\nu)}}=p$.
  Due to $m_1 = (n-2)p + (p-1)$, this is possible for at most $n-2$ different $\nu \in [1,n]$.
  Thus $k \ge 2$.

  We can also estimate $\lrcard{ \bigcup_{i=1}^4 \bigcup_{\chi \in B_i^{(\nu)}} \chi \langle g_i \rangle^\perp }$ in a different way:
  Assume for the purpose of showing \eqref{teetime} (the other cases are argued identically) that $n^{(\nu)} = \card{B_1^{(\nu)}} \ge \card{B_2^{(\nu)}} \ge \card{B_3^{(\nu)}} \ge \card{B_4^{(\nu)}}$.
  Each of the characters $\chi \in B_1^{(\nu)}$ contributes $\chi \langle g_1 \rangle^\perp$, and therefore exactly $p$ characters, to the union.
  Each of the characters $\chi \in B_2^{(\nu)}$ contributes at most $p - \card{B_1^{(\nu)}}$ characters, since $\card{\chi_1 \langle g_1 \rangle^\perp \cap \chi \langle g_2 \rangle^\perp} = 1$ for all $\chi_1 \in B_1^{(\nu)}$.
  Similarly, each of the characters $\chi \in B_3^{(\nu)}$ contributes at most $p - \max\{ \card{B_1^{(\nu)}} , \card{B_2^{(\nu)}} \} = p - \card{B_1^{(\nu)}}$ characters, since $\card{\chi_1 \langle g_1 \rangle^\perp \cap \chi \langle g_3 \rangle^\perp} = 1$ for all $\chi_1 \in B_1^{(\nu)}$ and $\card{\chi_2 \langle g_2 \rangle^\perp \cap \chi \langle g_3 \rangle^\perp} = 1$ for all $\chi_2 \in B_2^{(\nu)}$.
  Continuing this thought for $B_4^{(\nu)}$, we obtain
  \ber
  p^2 =\nn \lrcard{\bigcup_{i=1}^4 \bigcup_{\chi \in B_i^{(\nu)}} \chi \langle g_i \rangle^\perp} &
                \leq & p \card{B_1^{(\nu)}} + (p - \card{B_1^{(\nu)}})(\sum_{i=2}^4 \card{B_i^{(\nu)}}) \\
                & = &p n^{(\nu)} + (p - n^{(\nu)})(l^{(\nu)} - n^{(\nu)})\nn.
  \eer
  Therefore
  \be\label{teetime}
  (n^{(\nu)} - (l^{(\nu)} - p))(n^{(\nu)} - p) = p n^{(\nu)} + (p - n^{(\nu)})(l^{(\nu)} - n^{(\nu)}) - p^2 \ge 0.
  \ee
  Thus either $n^{(\nu)} \ge p$ (and therefore already $n^{(\nu)}=p$) or $n^{(\nu)} \le l^{(\nu)} - p$.

  For $\nu \in [1,k]$, we obtain $n^{(\nu)} \le l^{(\nu)} - p$.
  Due to $\card{B_4^{(\nu)}} \le m_4=2$, we also have $l^{(\nu)} = \sum_{i=1}^4 \card{B_i^{(\nu)}} \le 3n^{(\nu)} + 2$.
  Then
  \[
  3 l^{(\nu)} \geq 3n^{(\nu)} + 3p = 3 n^{(\nu)} + 2 + 3p - 2 \geq l^{(\nu)} + 3p - 2,
  \]
  and hence $l^{(\nu)} \geq \frac{3}{2} p - 1$ for all $\nu \in [1,k]$.
Because of $\sum_{i=1}^n l^{(\nu)} = \length{S} = pn + (p-1)$ and $l^{(\nu)} \geq n^{(\nu)} = p$ for all $\nu \in [k+1,n]$, we have $l^{(1)} + \ldots + l^{(k)} \le pk + (p-1)$.
For the remainder of the argument, we consider $\nu\in [1,k]$.

  Then, by the above, $\sum_{i=1, i\neq\nu}^{k} l^{(\nu)} \geq (k-1) (\frac{3}{2} p - 1)$, and hence
  \begin{equation}
    \label{eq:k-est}
  (k-1)\left(\frac{3}{2} p - 1\right) + l^{(\nu)} \leq pk + (p-1),
  \end{equation}
  which implies
  \[
  \begin{split}
  l^{(\nu)} & \leq pk + (p-1) - (k-1)\left(\frac{3}{2}p - 1\right)
          = pk + p - 1 - \frac{3}{2} kp  + k + \frac{3}{2} p - 1\\
          &= \frac{3}{2} p + (p - 2)  + k  - \frac{1}{2} pk
          = \frac{3}{2} p + (p - 2)  - \frac{k}{2}(p - 2)
  \end{split}
  \]
  Hence,  since $k \ge 2$, it follows that $l^{(\nu)} \le \lfloor \frac{3}{2}p \rfloor$. \footnote{Alternatively \eqref{eq:k-est}, together with $\l^{(\nu)} \ge \frac{3}{2}p - 1$, $p \ge 5$ and $k \le 2$, already implies $k = 2$, which yields the same estimate for $l^{(\nu)}$.}
  Together with $l^{(\nu)} \ge \lceil \frac{3}{2}p - 1 \rceil$, this implies $l^{(\nu)} = \frac{3}{2}p - \frac{1}{2}$.

  Since $\card{B_4^{(1)}} + \ldots + \card{B_4^{(k)}} \leq m_4 = 2$ and $k \ge 2$, there exists a $\nu \in [1,k]$ with $\card{B_4^{(\nu)}} \le 1$.
  Then
  \[
  \card{B_1^{(\nu)}},\ldots,\card{B_3^{(\nu)}} \leq n^{(\nu)} \leq l^{(\nu)} - p = \frac{1}{2} (p-1),
  \]
  $\card{B_4^{(\nu)}} \leq 1$ and $\sum_{i=1}^4 \card{B_i^{(\nu)}} = l^{(\nu)} = \frac{3}{2}(p-1)+1$.
  Therefore we must have $$\card{B_1^{(\nu)}}=\card{B_2^{(\nu)}}=\card{B_3^{(\nu)}} = n^{(\nu)}=\frac{1}{2}(p-1)$$ and $\card{B_4^{(\nu)}} = 1$.

  With the help of Lemma \ref{lemma:l-triple}, we show that this leads to a contradiction.
  Consider $T=g_1^{\frac{1}{2}(p-1)} g_2^{\frac{1}{2}(p-1)} g_3^{\frac{1}{2}(p-1)} \in \mathcal{F}(G)$.
  Then, by Lemma \ref{lemma:l-triple} (with $l=\frac{1}{2}(p-1)$ and $\chi=\varphi^{\sigma(\nu)}$),
  \[
  \lrcard{ \bigcup_{i=1}^3 \bigcup_{\chi' \in B_i^{(\nu)}} \chi' \langle g_i \rangle^\perp } < \frac{1}{2}(p-1)(2p+1).
  \]
  Thus, with $B_4^{(\nu)} = \{ \tau \}$,
  \[
  \begin{split}
  p^2 &= \lrcard{ \Big( \bigcup_{i=1}^3 \bigcup_{\chi' \in B_i^{(\nu)}} \chi' \langle g_i \rangle^\perp \Big) \cup \tau \langle g_4 \rangle^\perp }
      \le \lrcard{ \bigcup_{i=1}^3 \bigcup_{\chi' \in B_i^{(\nu)}} \chi' \langle g_i \rangle^\perp } + (p - n^{(\nu)}) \\
      &< \frac{1}{2}(p-1)(2p + 1) + \frac{1}{2}(p+1) = p^2,
  \end{split}
  \]
  a contradiction.
\end{proof}

\section{Acknowledgements}

I am indebted to Alfred Geroldinger for his constant feedback and help during the creation of this paper.
I would also like to thank David Grynkiewicz and Günter Lettl for their comments on preliminary versions of this paper.
In particular, G. Lettl suggested the use of the Cauchy-Davenport Theorem in Lemma \ref{lemma:l-triple}, which significantly shortened the proof, compared to an earlier version.


\providecommand{\bysame}{\leavevmode\hbox
to3em{\hrulefill}\thinspace}
\providecommand{\MR}{\relax\ifhmode\unskip\space\fi MR }
\providecommand{\MRhref}[2]{%
  \href{http://www.ams.org/mathscinet-getitem?mr=#1}{#2}
} \providecommand{\href}[2]{#2}

\end{document}